\documentclass[11pt,a4paper]{article}

\usepackage[utf8]{inputenc}
\usepackage[T1]{fontenc}
\usepackage{mathptmx}
\usepackage{amsmath,amsthm,amssymb}
\usepackage[margin=2.5cm]{geometry}
\usepackage{indentfirst}
\usepackage{hyperref}
\usepackage{setspace}
\usepackage{indentfirst}

\hypersetup{
    colorlinks=true,
    linkcolor=blue,
    citecolor=blue,
    urlcolor=blue
}

\newtheorem{definition}{Definition}[section]
\newtheorem{lemma}[definition]{Lemma}
\newtheorem{proposition}[definition]{Proposition}
\newtheorem{theorem}[definition]{Theorem}
\newtheorem{corollary}[definition]{Corollary}
\newtheorem{example}[definition]{Example}
\newtheorem{remark}[definition]{Remark}

\title{Properties of $f$-Correlated Fuzzy Numbers}
\author{Diogo Sampaio da Silva \thanks{diogo.sampaio@ufpel.edu.br} \and Roberto Antonio Cordeiro Prata}
\date{\today}

\begin{document}
\frenchspacing
\maketitle

\begin{abstract}
This paper investigates structural invariance and algebraic classification properties within the theory of interactive fuzzy numbers, focusing on $f$-correlated and more specifically, linearly correlated fuzzy numbers. Starting from the foundations of joint possibility distributions and the generalized extension principle, we analyze how continuous, strictly monotonic correlation functions $f$ act on $\alpha$-level sets and membership representations. We prove that $f$-correlation preserves the class $\mathcal{F}_C(\mathbb{R})$ of fuzzy numbers with continuous endpoints and transforms any $L$-$R$ fuzzy number into a generalized $L_f$-$R_f$ (or $R_f$-$L_f$) fuzzy number, with explicit formulas for the induced shape functions. Specializing to hyperbolic correlation ($f(x) = q/x + r$) illustrates how nonlinear interactivity preserves generalized $L$-$R$ regularity while deforming piecewise-linear triangular and trapezoidal flanks into rational curves. In contrast, when $f(x) = qx + r$ is affine, the reference shape functions $(L, R)$ are preserved up to orientation reversal, yielding complete set-theoretic characterizations of the $A$-linearly correlated spaces $\mathbb{R}_{\mathcal{F}(A)}$ for symmetric Gaussian, trapezoidal, and triangular generators, as well as necessary and sufficient shape-ratio conditions for linear correlation in the asymmetric setting.
\end{abstract}

\section{Introduction}

Introduced by Zadeh \cite{Zadeh1965}, fuzzy set theory generalizes classical set membership through continuous degrees in $[0, 1]$, providing a mathematical framework for modeling epistemic uncertainty and imprecise quantities. Extending real arithmetic to fuzzy numbers via Zadeh's extension principle implicitly assumes that the fuzzy quantities involved are \emph{non-interactive}, meaning that the realization of one variable imposes no constraint on the values assumed by another \cite{Barros2017Book,Fuller2004,Zadeh1975}. In mathematical modeling however, uncertain quantities frequently exhibit intrinsic functional or relational dependencies \cite{Cabral2015,Barros2017,Pedro2019}. Ignoring this interactivity might produce an artificial increase in uncertainty over time or between successive operations \cite{Barros2017,Pedro2019}. However, this is not merely a modeling approach: \cite{wasques2020interactive} demonstrated that the generalized fuzzy derivative for fuzzy numbers with continuous extremes is inherently interactive, so non-interactive approaches are inappropriate to fuzzy functions differentiable in this context.

To incorporate dependence structures without requiring complete specification of arbitrary multivariate possibility distributions, several classes of interactive fuzzy numbers have been proposed \cite{Carlsson2004,Coroianu2013a,Coroianu2013b,Fuller2004,Fuller2011}. Among these,  \emph{linearly correlated} fuzzy numbers restrict joint possibility distributions to affine lines $y = qx + r$ ($q \neq 0$) \cite{Carlsson2004,Barros2015,Barros2017}, while \emph{$f$-correlated} fuzzy numbers generalize this mechanism to arbitrary continuous, strictly monotonic injective maps $y = f(x)$, including hyperbolic interactions $f(x) = q/x + r$ \cite{Cabral2015}. These classes were used to develop an interactive differential calculus \cite{barros2017interactive}, which provides a suitable differential framework when fuzzy quantities are functionally related rather than independent.

Building on linear correlation, Esmi et al.\ \cite{Esmi2018} formalized the space $\mathbb{R}_{\mathcal{F}(A)}$ of all fuzzy numbers linearly correlated to a fixed generator $A$, establishing generator invariance $\mathbb{R}_{\mathcal{F}(A)} = \mathbb{R}_{\mathcal{F}(B)}$ for non-degenerate $B \in \mathbb{R}_{\mathcal{F}(A)} \setminus \mathbb{R}$ and developing a Fr\'{e}chet differential calculus on $\mathbb{R}_{\mathcal{F}(A)}$.

Despite these advances, two foundational questions remain essential for analytical and computational modeling: (i) which structural and parametric shape classes of fuzzy numbers are preserved under general nonlinear $f$-correlation versus linear correlation, and (ii) which explicit parametric subfamilies coincide with a given space $\mathbb{R}_{\mathcal{F}(A)}$? In this work, we synthesize the interactive extension principle \cite{Fuller2004,Cabral2015}, the affine operator framework of \cite{Barros2017,Esmi2018,Pedro2019}, and classical $L$-$R$ representations \cite{Dubois1980,Bede2013} to provide unified answers to both questions:
\begin{enumerate}
    \item We prove that any continuous, strictly monotonic $f$-correlation preserves the class $\mathcal{F}_C(\mathbb{R})$ of fuzzy numbers with continuous endpoints and maps strict or extended $L$-$R$ fuzzy numbers into generalized $L_f$-$R_f$ (for $f$ increasing) or $R_f$-$L_f$ (for $f$ decreasing) fuzzy numbers, deriving the exact transformed shape functions.
    \item Using hyperbolic correlation ($f(x) = q/x + r$) as a canonical nonlinear benchmark, we demonstrate explicitly that nonlinear $f$-correlation preserves generalized $L$-$R$ structure while breaking piecewise-linear trapezoidal and triangular geometry.
    \item Restricting to linear correlation ($f(x) = qx + r$), we establish exact set-theoretic characterizations of $\mathbb{R}_{\mathcal{F}(A)}$ for non-degenerate symmetric Gaussian, symmetric trapezoidal, and symmetric triangular generators, together with necessary and sufficient shape-invariance criteria for asymmetric Gaussian, trapezoidal, and triangular fuzzy numbers.
\end{enumerate}

\section{Theoretical Foundation}

Let $\mathcal{F}(\mathbb{R})$ denote the class of fuzzy subsets of $\mathbb{R}$ characterized by membership functions $\mu_A : \mathbb{R} \to [0, 1]$ \cite{Zadeh1965,Barros2017Book}. For $\alpha \in (0, 1]$, the $\alpha$-level set of $A \in \mathcal{F}(\mathbb{R})$ is defined as
\begin{equation}
[A]_\alpha = \{x \in \mathbb{R} : \mu_A(x) \ge \alpha\},
\end{equation}
and $[A]_0 = \overline{\{x \in \mathbb{R} : \mu_A(x) > 0\}}$ denotes the closure of its support \cite{Barros2017Book,Pedro2019}. A fundamental property of $\alpha$-levels is that they uniquely determine fuzzy subsets.

\begin{lemma}[Level-Set Uniqueness \cite{Barros2017Book,Esmi2018}]
\label{lem:uniqueness}
Let $A, B \in \mathcal{F}(\mathbb{R})$. Then $A = B$ if and only if $[A]_\alpha = [B]_\alpha$ for all $\alpha \in (0, 1]$ (equivalently, for all $\alpha \in [0, 1]$).
\end{lemma}

\begin{proof}
Clearly $A = B$ implies $[A]_\alpha = [B]_\alpha$ for all $\alpha \in (0, 1]$. Conversely, suppose $A \neq B$, so there exists $x_0 \in \mathbb{R}$ with $\mu_A(x_0) \neq \mu_B(x_0)$. Without loss of generality, assume $\mu_A(x_0) > \mu_B(x_0)$ and set $\alpha_0 = \mu_A(x_0) \in (0, 1]$. Then $x_0 \in [A]_{\alpha_0}$ whereas $x_0 \notin [B]_{\alpha_0}$, contradicting $[A]_{\alpha_0} = [B]_{\alpha_0}$.
\end{proof}

\begin{definition}[Fuzzy Number \cite{Barros2017Book,Esmi2018}]
\label{def:fuzzy_number}
A fuzzy subset $A \in \mathcal{F}(\mathbb{R})$ is called a \emph{fuzzy number} (denoted $A \in \mathbb{R}_{\mathcal{F}}$) if $[A]_\alpha := [a_\alpha^-, a_\alpha^+] \subset \mathbb{R}$ is a non-empty bounded closed interval for every $\alpha \in [0, 1]$ (or for every $\alpha \in (0, 1]$ when unbounded supports such as classical Gaussian membership functions are admitted). Every crisp real number $r \in \mathbb{R}$ is embedded in $\mathbb{R}_{\mathcal{F}}$ via its characteristic function $\chi_{\{r\}}$, whose $\alpha$-levels satisfy $[\chi_{\{r\}}]_\alpha = \{r\}$ for all $\alpha \in [0, 1]$.
\end{definition}

To model dependencies among uncertain variables, we recall the notions of possibility distributions and interactive extensions \cite{Fuller2004,Fuller2011}.

\begin{definition}[Possibility and Joint Possibility Distributions \cite{Fuller2004,Fuller2011}]
\label{def:joint_poss}
A function $\mu_C : \mathbb{R}^n \to [0, 1]$ is a \emph{joint possibility distribution} of the fuzzy numbers $A_1, \dots, A_n \in \mathbb{R}_{\mathcal{F}}$ if $\sup_{x \in \mathbb{R}^n} \mu_C(x) = 1$ and each $\mu_{A_i}$ is recovered as its marginal distribution:
\begin{equation}
\mu_{A_i}(x_i) = \sup_{x_j \in \mathbb{R},\; j \neq i} \mu_C(x_1, \dots, x_n), \qquad i = 1, \dots, n.
\end{equation}
The fuzzy numbers $A_1, \dots, A_n$ are said to be \emph{non-interactive} if $\mu_C(x_1, \dots, x_n) = \min\{\mu_{A_1}(x_1), \dots, \mu_{A_n}(x_n)\}$ for all $(x_1, \dots, x_n) \in \mathbb{R}^n$; otherwise, they are called \emph{interactive}.
\end{definition}

\begin{definition}[Extension Principle with Joint Possibility Distribution \cite{Fuller2004}]
\label{def:ext_principle}
Let $C$ be a joint possibility distribution of $A_1, \dots, A_n \in \mathbb{R}_{\mathcal{F}}$, and let $g : \mathbb{R}^n \to \mathbb{R}$ be a continuous function. The interactive extension $g_C(A_1, \dots, A_n) \in \mathcal{F}(\mathbb{R})$ is defined by
\begin{equation}
\mu_{g_C(A_1, \dots, A_n)}(y) = \sup_{y = g(x_1, \dots, x_n)} \mu_C(x_1, \dots, x_n), \qquad y \in \mathbb{R},
\end{equation}
with the convention $\sup \emptyset = 0$.
\end{definition}

This extension principle is the one used by the interactive derivative \cite{barros2017interactive}, and it is the differential counterpart of the level-wise representation developed in Section 3

\begin{lemma}[\cite{Fuller2004}]
\label{lem:nguyen_joint}
Let $A_1, \dots, A_n \in \mathbb{R}_{\mathcal{F}}$ have compact joint $\alpha$-levels $[C]_\alpha \subset \mathbb{R}^n$, and let $g : \mathbb{R}^n \to \mathbb{R}$ be continuous. Then, for every $\alpha \in [0, 1]$,
\begin{equation}
[g_C(A_1, \dots, A_n)]_\alpha = g([C]_\alpha).
\end{equation}
\end{lemma}

\begin{proof}
For any \( \alpha \in [0, 1] \) we obtain:
\begin{align}
\left [ g_{C} (A_1, \ldots, A_n) \right ]^{\alpha} &= \overline{\left\{ y \in \mathbb{R} ; \sup_{y=g(x)}  C(x) > \alpha \right\}}\\ 
 &= \overline{\left\{ y \in \mathbb{R} ; \exists x \in \mathbb{R} ; y = g(x), C(x) > \alpha \right\}}\\ 
 &= \overline{\left\{ g(x) \in \mathbb{R} ; C(x) > \alpha \right\}}\\ 
 &= g \left ( \overline{\left\{ x \in \mathbb{R}^{n} ; C(x) > \alpha \right\}} \right ) \\
 &= g([C]^{\alpha })
\end{align}
This concludes our proof.
\end{proof}

When two fuzzy numbers are functionally coupled through a continuous, monotonic injective function, their joint possibility distribution concentrates on the graph of that function \cite{Cabral2015}.

\begin{definition}[$f$-Correlated Fuzzy Numbers \cite{Cabral2015}]
\label{def:f_corr}
Let $X, Y \subseteq \mathbb{R}$ be intervals and let $f : X \to Y$ be a continuous, strictly monotonic (hence injective) function. Two fuzzy numbers $A, B \in \mathbb{R}_{\mathcal{F}}$ are said to be \emph{$f$-correlated} if $[A]_\alpha \subset X$ for all $\alpha \in (0, 1]$ (with $[A]_0 \subset X$ when $[A]_0$ is compact) and their joint possibility distribution $J$ is given by
\begin{equation}
\mu_J(x, y) = \mu_A(x)\chi_{\{y = f(x)\}}(x, y) = \mu_B(y)\chi_{\{y = f(x)\}}(x, y), \qquad (x, y) \in \mathbb{R}^2,
\end{equation}
where $\chi_{\{y = f(x)\}}$ denotes the characteristic function of the graph $\{(x, f(x)) : x \in X\}$. In this case, $[J]_\alpha = \{(x, f(x)) : x \in [A]_\alpha\}$ and $[B]_\alpha = f([A]_\alpha)$ for all $\alpha \in (0, 1]$ (and $\alpha = 0$ when $[A]_0 \subset X$).
\end{definition}

\begin{example}[Hyperbolically Interactive Fuzzy Numbers \cite{Cabral2015}]
\label{ex:hyperbolic}
Let $q, r \in \mathbb{R}$ with $q \neq 0$, and consider the hyperbolic map $f(x) = \frac{q}{x} + r$ defined on either $X = (0, +\infty)$ or $X = (-\infty, 0)$. Two fuzzy numbers $A$ and $B$ with $[A]_0 \subset X$ are \emph{hyperbolically interactive} if they are $f$-correlated with joint possibility distribution $\mu_J(x, y) = \mu_A(x)\chi_{\{y = q/x + r\}}(x, y)$, so that
\begin{equation}
[B]_\alpha = \left\{\frac{q}{x} + r : x \in [A]_\alpha\right\}, \qquad \alpha \in [0, 1].
\end{equation}
\end{example}

\begin{example}[Linearly Correlated Fuzzy Numbers and $\mathbb{R}_{\mathcal{F}(A)}$ \cite{Carlsson2004,Barros2017,Esmi2018}]
\label{def:lin_corr}
Two fuzzy numbers $A, B \in \mathbb{R}_{\mathcal{F}}$ are said to be \emph{linearly correlated} (or \emph{completely correlated} when $q \neq 0$) if there exist scalars $q, r \in \mathbb{R}$ such that
\begin{equation}
[B]_\alpha = q[A]_\alpha + r
\end{equation}
for all $\alpha \in (0, 1]$ (and $\alpha \in [0, 1]$ when $[A]_0$ is compact). Following Esmi et al.\ \cite{Esmi2018}, for a fixed generator $A \in \mathbb{R}_{\mathcal{F}}$, define the operator $\Psi_A : \mathbb{R}^2 \to \mathbb{R}_{\mathcal{F}}$ by $[\Psi_A(q, r)]_\alpha = q[A]_\alpha + r$. Its range
\begin{equation}
\mathbb{R}_{\mathcal{F}(A)} = \{\Psi_A(q, r) : (q, r) \in \mathbb{R}^2\}
\end{equation}
is the space of $A$-linearly correlated fuzzy numbers. Notably, $\mathbb{R} \subseteq \mathbb{R}_{\mathcal{F}(A)}$ since $\chi_{\{r\}} = \Psi_A(0, r)$ for every $r \in \mathbb{R}$ \cite{Esmi2018,Pedro2019}.
\end{example}

We also recall two results from  \cite{Esmi2018} concerning generator invariance and symmetry.

\begin{proposition}[Generator Invariance \cite{Esmi2018}]
\label{prop:invariance}
Let $B = \Psi_A(q, r) \in \mathbb{R}_{\mathcal{F}}$ for some $A \in \mathbb{R}_{\mathcal{F}}$ and $(q, r) \in \mathbb{R}^2$ with $q \neq 0$. Then $\mathbb{R}_{\mathcal{F}(A)} = \mathbb{R}_{\mathcal{F}(B)}$.
\end{proposition}

\begin{lemma}[Level-wise Symmetry and Affine Preservation \cite{Esmi2018}]
\label{lem:symmetry}
A fuzzy number $A \in \mathbb{R}_{\mathcal{F}}$ is symmetric with respect to $x \in \mathbb{R}$ (i.e., $\mu_A(x - y) = \mu_A(x + y)$ for all $y \in \mathbb{R}$) if and only if $a_\alpha^+ = 2x - a_\alpha^-$ (equivalently, $x = \frac{1}{2}(a_\alpha^- + a_\alpha^+)$) for all $\alpha \in (0, 1]$. Moreover, if $A \in \mathbb{R}_{\mathcal{F}} \setminus \mathbb{R}$ is symmetric with respect to $x$, then for any $(q, r) \in \mathbb{R}^2$, $\Psi_A(q, r) = \Psi_A(-q, 2qx + r)$ is symmetric with respect to $qx + r$.
\end{lemma}

\section{Structural Preservation under $f$-Correlation}

We begin by characterizing the action of a general $f$-correlation on compact $\alpha$-level intervals and proving that it preserves both endpoint continuity and generalized $L$-$R$ membership representations.

\begin{lemma}[Image of Compact Intervals]
\label{lem:interval_image}
Let $X \subseteq \mathbb{R}$ be an interval, $[a, b] \subset X$, and $f : X \to \mathbb{R}$ be continuous and strictly monotonic on $[a, b]$. Then:
\begin{enumerate}
    \item If $f$ is strictly increasing, $f([a, b]) = [f(a), f(b)]$.
    \item If $f$ is strictly decreasing, $f([a, b]) = [f(b), f(a)]$.
\end{enumerate}
\end{lemma}

\begin{proof}
Suppose $f$ is strictly increasing on $[a, b]$. For every $x \in [a, b]$, $a \le x \le b$ implies $f(a) \le f(x) \le f(b)$, so $f([a, b]) \subseteq [f(a), f(b)]$. Conversely, for any $y \in [f(a), f(b)]$, continuity of $f$ on $[a, b]$ and the Intermediate Value Theorem guarantee the existence of $c \in [a, b]$ such that $f(c) = y$, proving $[f(a), f(b)] \subseteq f([a, b])$. The strictly decreasing case follows identically upon reversing the endpoint inequalities $f(b) \le f(x) \le f(a)$.
\end{proof}

\begin{theorem}[Level-Set Representation of $f$-Correlated Fuzzy Numbers]
\label{thm:f_corr_levels}
Let $A \in \mathbb{R}_{\mathcal{F}}$ have $\alpha$-levels $[A]_\alpha = [a_\alpha^-, a_\alpha^+]$, and let $B \in \mathbb{R}_{\mathcal{F}}$ be $f$-correlated to $A$ via a continuous, strictly monotonic function $f : X \to Y$. Then, for all $\alpha \in (0, 1]$ (and $\alpha = 0$ when $[A]_0 \subset X$):
\begin{enumerate}
    \item If $f$ is strictly increasing, $[B]_\alpha = [f(a_\alpha^-), f(a_\alpha^+)]$.
    \item If $f$ is strictly decreasing, $[B]_\alpha = [f(a_\alpha^+), f(a_\alpha^-)]$.
\end{enumerate}
\end{theorem}

\begin{proof}
Since $[B]_\alpha = f([A]_\alpha) = f([a_\alpha^-, a_\alpha^+])$ by Definition~\ref{def:f_corr} and $[a_\alpha^-, a_\alpha^+] \subset X$ is a compact interval, the result follows directly from Lemma~\ref{lem:interval_image}.
\end{proof}

We now recall two important classes of fuzzy numbers \cite{Dubois1980,Bede2013}:
\begin{enumerate}
    \item \emph{Fuzzy numbers with continuous endpoints ($\mathcal{F}_C(\mathbb{R})$):} A fuzzy number $A \in \mathbb{R}_{\mathcal{F}}$ with $[A]_\alpha = [a_\alpha^-, a_\alpha^+]$ belongs to $\mathcal{F}_C(\mathbb{R})$ if the endpoint mappings $\alpha \mapsto a_\alpha^-$ and $\alpha \mapsto a_\alpha^+$ are continuous on $[0, 1]$ (or on $(0, 1]$ in the unbounded-support setting) \cite{Bede2013}.
    \item \emph{$L$-$R$ fuzzy numbers (strict and extended):} Given continuous, non-increasing reference shape functions $L, R : [0, \infty) \to [0, 1]$ with $L(0) = R(0) = 1$, a \emph{strict $L$-$R$ fuzzy number} $A = (m, \alpha_L, \beta_R)_{LR}$ with modal value $m \in \mathbb{R}$ and spreads $\alpha_L, \beta_R > 0$ has membership function \cite{Dubois1980,Bede2013}
    \begin{equation}
    \mu_A(x) = \begin{cases}
    L\!\left(\dfrac{m - x}{\alpha_L}\right), & x \le m, \\[10pt]
    R\!\left(\dfrac{x - m}{\beta_R}\right), & x > m.
    \end{cases}
    \end{equation}
    Under the \emph{extended $L$-$R$ convention}, zero spreads ($\alpha_L, \beta_R \ge 0$) and flat cores $[m_1, m_2]$ ($m_1 \le m_2$) are permitted by setting the corresponding flank to $0$ outside $[m_1, m_2]$ whenever a spread vanishes, thereby encompassing crisp real numbers $(r, 0, 0)_{LR} = \chi_{\{r\}}$ and crisp intervals $(a, b, 0, 0)_{LR} = \chi_{[a, b]}$ \cite{Dubois1980,Bede2013}. When $L$ and $R$ are supported on $[0, 1]$ with $L(1) = R(1) = 0$, the support of $A = (m_1, m_2, \alpha_L, \beta_R)_{LR}$ is $[A]_0 = [m_1 - \alpha_L, m_2 + \beta_R]$.
\end{enumerate}

\begin{theorem}[Preservation of Continuous Endpoints and Generalized $L$-$R$ Type under $f$-Correlation]
\label{thm:f_corr_preservation}
Let $A, B \in \mathbb{R}_{\mathcal{F}}$ be $f$-correlated via a continuous, strictly monotonic function $f : X \to Y$. Then:
\begin{enumerate}
    \item \emph{(Continuous Endpoints)} If $A \in \mathcal{F}_C(\mathbb{R})$, then $B \in \mathcal{F}_C(\mathbb{R})$.
    \item \emph{(Generalized $L$-$R$ Type)} Suppose $A = (m_1, m_2, \alpha_L, \beta_R)_{LR}$ is a non-degenerate ($\alpha_L, \beta_R > 0$) $L$-$R$ fuzzy number with continuous non-increasing reference functions $L, R : [0, 1] \to [0, 1]$ satisfying $L(0) = R(0) = 1$, $L(1) = R(1) = 0$, and support $[A]_0 = [m_1 - \alpha_L, m_2 + \beta_R] \subset X$. Then $B$ is also an $L$-$R$ fuzzy number:
    \begin{itemize}
        \item If $f$ is strictly increasing, then $B = (\tilde{m}_1, \tilde{m}_2, \tilde{\alpha}_L, \tilde{\beta}_R)_{L_f R_f}$ with $\tilde{m}_1 = f(m_1)$, $\tilde{m}_2 = f(m_2)$, $\tilde{\alpha}_L = f(m_1) - f(m_1 - \alpha_L) > 0$, $\tilde{\beta}_R = f(m_2 + \beta_R) - f(m_2) > 0$, and continuous non-increasing reference functions $L_f, R_f : [0, 1] \to [0, 1]$ given by
        \begin{align}
        L_f(u) &= L\!\left(\frac{m_1 - f^{-1}\!\bigl(f(m_1) - u\tilde{\alpha}_L\bigr)}{\alpha_L}\right), \label{eq:Lf_inc} \\
        R_f(u) &= R\!\left(\frac{f^{-1}\!\bigl(f(m_2) + u\tilde{\beta}_R\bigr) - m_2}{\beta_R}\right), \qquad u \in [0, 1]. \label{eq:Rf_inc}
        \end{align}
        \item If $f$ is strictly decreasing, then $B = (\hat{m}_1, \hat{m}_2, \hat{\alpha}_L, \hat{\beta}_R)_{\tilde{L}_f \tilde{R}_f}$ with $\hat{m}_1 = f(m_2)$, $\hat{m}_2 = f(m_1)$, $\hat{\alpha}_L = f(m_2) - f(m_2 + \beta_R) > 0$, $\hat{\beta}_R = f(m_1 - \alpha_L) - f(m_1) > 0$, and continuous non-increasing reference functions $\tilde{L}_f, \tilde{R}_f : [0, 1] \to [0, 1]$ given by
        \begin{align}
        \tilde{L}_f(u) &= R\!\left(\frac{f^{-1}\!\bigl(f(m_2) - u\hat{\alpha}_L\bigr) - m_2}{\beta_R}\right), \label{eq:Lf_dec} \\
        \tilde{R}_f(u) &= L\!\left(\frac{m_1 - f^{-1}\!\bigl(f(m_1) + u\hat{\beta}_R\bigr)}{\alpha_L}\right), \qquad u \in [0, 1]. \label{eq:Rf_dec}
        \end{align}
    \end{itemize}
\end{enumerate}
\end{theorem}

\begin{proof}
\begin{enumerate}
    \item By Theorem~\ref{thm:f_corr_levels}, the endpoints $[b_\alpha^-, b_\alpha^+]$ of $B$ are given by $[f(a_\alpha^-), f(a_\alpha^+)]$ when $f$ is strictly increasing and $[f(a_\alpha^+), f(a_\alpha^-)]$ when $f$ is strictly decreasing. Since $f$ is continuous and $\alpha \mapsto a_\alpha^-$, $\alpha \mapsto a_\alpha^+$ are continuous, their compositions $\alpha \mapsto b_\alpha^-$ and $\alpha \mapsto b_\alpha^+$ are continuous, so $B \in \mathcal{F}_C(\mathbb{R})$.
    \item Because $f : X \to f(X)$ is continuous and strictly monotonic on the interval $X$, $f$ is a homeomorphism onto the interval $f(X)$ with continuous, strictly monotonic inverse $f^{-1} : f(X) \to X$. By Definition~\ref{def:f_corr}, $\mu_B(y) = \mu_A(f^{-1}(y))$ for $y \in f([A]_0)$ and $\mu_B(y) = 0$ otherwise.
    
    Assume $f$ is strictly increasing. Then $f([A]_0) = [f(m_1 - \alpha_L), f(m_2 + \beta_R)] = [\tilde{m}_1 - \tilde{\alpha}_L, \tilde{m}_2 + \tilde{\beta}_R]$ with $\tilde{\alpha}_L, \tilde{\beta}_R > 0$. For $y \in [\tilde{m}_1, \tilde{m}_2]$, we have $f^{-1}(y) \in [m_1, m_2]$ and thus $\mu_B(y) = 1$. For $y \in [\tilde{m}_1 - \tilde{\alpha}_L, \tilde{m}_1]$, setting $u = (\tilde{m}_1 - y)/\tilde{\alpha}_L \in [0, 1]$ gives $y = f(m_1) - u\tilde{\alpha}_L$ and $f^{-1}(y) \in [m_1 - \alpha_L, m_1]$, so
    \begin{equation}
    \mu_B(y) = \mu_A(f^{-1}(y)) = L\!\left(\frac{m_1 - f^{-1}(f(m_1) - u\tilde{\alpha}_L)}{\alpha_L}\right) = L_f\!\left(\frac{\tilde{m}_1 - y}{\tilde{\alpha}_L}\right).
    \end{equation}
    Moreover, $L_f(0) = L(0) = 1$, $L_f(1) = L(1) = 0$, and as $u$ increases on $[0, 1]$, $f(m_1) - u\tilde{\alpha}_L$ strictly decreases, $f^{-1}(f(m_1) - u\tilde{\alpha}_L)$ strictly decreases, and hence $L_f(u)$ is continuous and non-increasing. An identical argument on $[\tilde{m}_2, \tilde{m}_2 + \tilde{\beta}_R]$ yields $\mu_B(y) = R_f((y - \tilde{m}_2)/\tilde{\beta}_R)$ with $R_f(0) = 1$, $R_f(1) = 0$, and $R_f$ continuous and non-increasing.
    
    When $f$ is strictly decreasing, $f([A]_0) = [f(m_2 + \beta_R), f(m_1 - \alpha_L)] = [\hat{m}_1 - \hat{\alpha}_L, \hat{m}_2 + \hat{\beta}_R]$ with $\hat{\alpha}_L, \hat{\beta}_R > 0$. The left flank of $B$ corresponds to $f([m_2, m_2 + \beta_R])$ and is governed by $R$, while the right flank corresponds to $f([m_1 - \alpha_L, m_1])$ and is governed by $L$. Substituting $y = f(m_2) - u\hat{\alpha}_L$ on $[\hat{m}_1 - \hat{\alpha}_L, \hat{m}_1]$ and $y = f(m_1) + u\hat{\beta}_R$ on $[\hat{m}_2, \hat{m}_2 + \hat{\beta}_R]$ yields \eqref{eq:Lf_dec}--\eqref{eq:Rf_dec} with $\tilde{L}_f(0) = \tilde{R}_f(0) = 1$, $\tilde{L}_f(1) = \tilde{R}_f(1) = 0$, and both functions continuous and non-increasing.
\end{enumerate}
\end{proof}

When $f$ is affine, the reference shape functions $(L, R)$ are invariant without distortion, recovering and extending the linear correlation results of \cite{Barros2017,Esmi2018,Pedro2019}.

\begin{corollary}[Preservation under Linear Correlation \cite{Esmi2018,Pedro2019}]
\label{cor:linear_preservation}
Let $A \in \mathbb{R}_{\mathcal{F}}$ and $B = \Psi_A(q, r) \in \mathbb{R}_{\mathcal{F}(A)}$ for $(q, r) \in \mathbb{R}^2$. Then:
\begin{enumerate}
    \item If $A \in \mathcal{F}_C(\mathbb{R})$, then $B \in \mathcal{F}_C(\mathbb{R})$ for all $(q, r) \in \mathbb{R}^2$.
    \item If $A = (m_A, \alpha_L, \beta_R)_{LR}$ is a strict $L$-$R$ fuzzy number ($\alpha_L, \beta_R > 0$) and $q \neq 0$, then
    \begin{equation}
    B = \begin{cases}
    (q m_A + r,\; q\alpha_L,\; q\beta_R)_{LR}, & q > 0, \\
    (q m_A + r,\; |q|\beta_R,\; |q|\alpha_L)_{RL}, & q < 0.
    \end{cases}
    \end{equation}
    In particular, if $L = R$, then $B = (q m_A + r, |q|\alpha_L, |q|\beta_R)_{LL}$ preserves $(L, L)$ for all $q \neq 0$.
    \item If $A = (m_1, m_2, \alpha_L, \beta_R)_{LR}$ is an extended $L$-$R$ fuzzy number ($\alpha_L, \beta_R \ge 0$), then $B \in \mathbb{R}_{\mathcal{F}(A)}$ is an extended $L$-$R$ fuzzy number for all $(q, r) \in \mathbb{R}^2$, including $q = 0$ where $B = (r, r, 0, 0)_{LR} = \chi_{\{r\}}$.
\end{enumerate}
\end{corollary}

\begin{proof}
For $q = 0$, $[B]_\alpha = \{r\}$ has constant (hence continuous) endpoints and equals $(r, r, 0, 0)_{LR} = \chi_{\{r\}}$.

For $q \neq 0$, apply Theorem~\ref{thm:f_corr_preservation} (or direct affine substitution $\mu_B(x) = \mu_A((x - r)/q)$ \cite{Pedro2019} when $[A]_0$ is unbounded) with $f(x) = qx + r$ and $f^{-1}(y) = (y - r)/q$. 

For $q > 0$, $\tilde{\alpha}_L = q\alpha_L$ and $f^{-1}(f(m_1) - u q\alpha_L) = m_1 - u\alpha_L$, so \eqref{eq:Lf_inc}--\eqref{eq:Rf_inc} simplify to $L_f(u) = L(u)$ and $R_f(u) = R(u)$. 

For $q < 0$, $\hat{\alpha}_L = |q|\beta_R$ and $\hat{\beta}_R = |q|\alpha_L$, so \eqref{eq:Lf_dec}--\eqref{eq:Rf_dec} simplify to $\tilde{L}_f(u) = R(u)$ and $\tilde{R}_f(u) = L(u)$.
\end{proof}

\begin{example}[Hyperbolic Correlation of a Triangular Fuzzy Number]
\label{ex:hyperbolic_tri}
To illustrate the contrast between nonlinear $f$-correlation and linear correlation, let $A = (a; b; d)$ be a non-degenerate triangular fuzzy number with $0 < a < b < d$, represented as $A = (b, b, \alpha_L, \beta_R)_{LL}$ with $\alpha_L = b - a > 0$, $\beta_R = d - b > 0$, and $L(u) = R(u) = 1 - u$ for $u \in [0, 1]$. Suppose $B$ is hyperbolically correlated to $A$ via $f(x) = \frac{q}{x} + r$ on $X = (0, +\infty)$ with $q > 0$ (Example~\ref{ex:hyperbolic}). Since $f$ is continuous and strictly decreasing on $(0, +\infty)$ with $f^{-1}(y) = \frac{q}{y - r}$, Theorem~\ref{thm:f_corr_preservation}(2) implies that $B$ is a strict $\tilde{L}_f$-$\tilde{R}_f$ fuzzy number with modal peak $\hat{m} = \frac{q}{b} + r$, support $\left[\frac{q}{d} + r,\; \frac{q}{a} + r\right]$, and spreads
\begin{equation}
\hat{\alpha}_L = \frac{q}{b} - \frac{q}{d} = \frac{q(d - b)}{bd}, \qquad \hat{\beta}_R = \frac{q}{a} - \frac{q}{b} = \frac{q(b - a)}{ab}.
\end{equation}
Evaluating \eqref{eq:Lf_dec}--\eqref{eq:Rf_dec} with $f^{-1}(f(b) - u\hat{\alpha}_L) = \frac{bd}{d - (d - b)u}$ and $f^{-1}(f(b) + u\hat{\beta}_R) = \frac{ab}{a + (b - a)u}$ yields the induced reference shape functions
\begin{equation}
\tilde{L}_f(u) = \frac{d(1 - u)}{d - (d - b)u}, \qquad \tilde{R}_f(u) = \frac{a(1 - u)}{a + (b - a)u}, \qquad u \in [0, 1].
\end{equation}
Both $\tilde{L}_f$ and $\tilde{R}_f$ are continuous and strictly decreasing on $[0, 1]$ with $\tilde{L}_f(0) = \tilde{R}_f(0) = 1$ and $\tilde{L}_f(1) = \tilde{R}_f(1) = 0$, confirming $B \in \mathcal{F}_C(\mathbb{R})$ and generalized $L$-$R$ membership. However, because $\tilde{L}_f$ and $\tilde{R}_f$ are non-affine rational functions, $B$ is \emph{not} a triangular fuzzy number---demonstrating that piecewise-linear triangular and trapezoidal shapes are invariant under linear correlation (Corollary~\ref{cor:linear_preservation}) but not under general nonlinear $f$-correlation.
\end{example}

\begin{remark}
\label{rem:converse_failure}
While linear correlation preserves both $\mathcal{F}_C(\mathbb{R})$ and the reference pair $(L, R)$ up to transposition (Corollary~\ref{cor:linear_preservation}), two symmetric $L$-$R$ fuzzy numbers generated by distinct reference shape functions (for instance, triangular $L_1(u) = \max(0, 1-u)$ and Gaussian $L_2(u) = e^{-u^2/2}$) cannot be linearly correlated because $L_2^{-1}(\alpha)/L_1^{-1}(\alpha)$ is non-constant on $(0, 1)$.
\end{remark}

\section{Characterization of $\mathbb{R}_{\mathcal{F}(A)}$ for Symmetric Generators}

We now determine the exact set-theoretic structure of $\mathbb{R}_{\mathcal{F}(A)}$ when the generator $A$ belongs to the classical symmetric Gaussian, trapezoidal, or triangular families. The following characterizations should be read alongside the approximation framework of \cite{bello2026}: the characterizations say exactly which LR-type numbers lie in $\mathbb{R}_{\mathcal{F}(A)}$, while the approximation framework says how closely an arbitrary fuzzy number can be approximated by one of them.

\begin{definition}[Symmetric Gaussian Fuzzy Number]
\label{def:sym_gauss}
A fuzzy number $A$ is called a \emph{non-degenerate symmetric Gaussian fuzzy number} with center $m_A \in \mathbb{R}$ and spread $\sigma_A > 0$ if its membership function is given by
\begin{equation}
\mu_A(x) = \exp\!\left(-\frac{(x - m_A)^2}{2\sigma_A^2}\right), \qquad x \in \mathbb{R}.
\end{equation}
We denote by $\mathcal{G}_{\mathrm{sym}}$ the family of all non-degenerate symmetric Gaussian fuzzy numbers.
\end{definition}

\begin{theorem}[Characterization of $\mathbb{R}_{\mathcal{F}(A)}$ for Symmetric Gaussians]
\label{prop:gauss_char}
Let $A \in \mathcal{G}_{\mathrm{sym}}$ be arbitrary. Then:
\begin{enumerate}
    \item Every $B \in \mathcal{G}_{\mathrm{sym}}$ is linearly correlated to $A$;
    \item $\mathbb{R}_{\mathcal{F}(A)} = \mathcal{G}_{\mathrm{sym}} \cup \mathbb{R}$, and consequently $\mathbb{R}_{\mathcal{F}(A)} = \mathbb{R}_{\mathcal{F}(B)}$ for all $A, B \in \mathcal{G}_{\mathrm{sym}}$.
\end{enumerate}
\end{theorem}

\begin{proof}
Let $A \in \mathcal{G}_{\mathrm{sym}}$ have center $m_A \in \mathbb{R}$ and spread $\sigma_A > 0$. For any $\alpha \in (0, 1]$, solving $\mu_A(x) \ge \alpha$ yields
\begin{equation}
\exp\!\left(-\frac{(x - m_A)^2}{2\sigma_A^2}\right) \ge \alpha \iff (x - m_A)^2 \le -2\sigma_A^2 \ln\alpha \iff |x - m_A| \le \sigma_A \sqrt{-2\ln\alpha},
\end{equation}
so the $\alpha$-levels of $A$ are
\begin{equation}
[A]_\alpha = \left[m_A - \sigma_A\sqrt{-2\ln\alpha},\; m_A + \sigma_A\sqrt{-2\ln\alpha}\right], \qquad \alpha \in (0, 1].
\end{equation}

To establish $\mathcal{G}_{\mathrm{sym}} \cup \mathbb{R} \subseteq \mathbb{R}_{\mathcal{F}(A)}$, let $B \in \mathcal{G}_{\mathrm{sym}}$ have center $m_B \in \mathbb{R}$ and spread $\sigma_B > 0$, with
\begin{equation}
[B]_\alpha = \left[m_B - \sigma_B\sqrt{-2\ln\alpha},\; m_B + \sigma_B\sqrt{-2\ln\alpha}\right], \qquad \alpha \in (0, 1].
\end{equation}
Setting $q = \sigma_B / \sigma_A > 0$ and $r = m_B - (\sigma_B / \sigma_A)m_A$, interval scaling preserves endpoint ordering \cite{Barros2017,Esmi2018} and gives, for all $\alpha \in (0, 1]$,
\begin{align}
[\Psi_A(q, r)]_\alpha &= q[A]_\alpha + r \nonumber \\
&= \left[\frac{\sigma_B}{\sigma_A}m_A + m_B - \frac{\sigma_B}{\sigma_A}m_A - \sigma_B\sqrt{-2\ln\alpha},\; \frac{\sigma_B}{\sigma_A}m_A + m_B - \frac{\sigma_B}{\sigma_A}m_A + \sigma_B\sqrt{-2\ln\alpha}\right] \nonumber \\
&= \left[m_B - \sigma_B\sqrt{-2\ln\alpha},\; m_B + \sigma_B\sqrt{-2\ln\alpha}\right] = [B]_\alpha.
\end{align}
Thus $B = \Psi_A(q, r) \in \mathbb{R}_{\mathcal{F}(A)}$, proving item~(1). Since $[\chi_{\{r\}}]_\alpha = \{r\} = [\Psi_A(0, r)]_\alpha$ for every $r \in \mathbb{R}$ \cite{Esmi2018,Pedro2019}, we obtain $\mathcal{G}_{\mathrm{sym}} \cup \mathbb{R} \subseteq \mathbb{R}_{\mathcal{F}(A)}$.

Conversely, let $C = \Psi_A(q, r) \in \mathbb{R}_{\mathcal{F}(A)}$ for arbitrary $(q, r) \in \mathbb{R}^2$. Because $A$ is a strict symmetric $L$-$L$ fuzzy number with $L(u) = e^{-u^2/2}$ and $\alpha_L = \beta_R = \sigma_A$, Corollary~\ref{cor:linear_preservation} implies that if $q = 0$, then $C = \chi_{\{r\}} \in \mathbb{R}$, whereas if $q \neq 0$, then $C$ is a strict symmetric $L$-$L$ fuzzy number with reference shape $L(u) = e^{-u^2/2}$, center $m_C = q m_A + r$, and spread $\sigma_C = |q|\sigma_A > 0$, i.e., $C \in \mathcal{G}_{\mathrm{sym}}$ (consistent with Lemma~\ref{lem:symmetry}). Hence $\mathbb{R}_{\mathcal{F}(A)} = \mathcal{G}_{\mathrm{sym}} \cup \mathbb{R}$, and Proposition~\ref{prop:invariance} ensures $\mathbb{R}_{\mathcal{F}(A)} = \mathbb{R}_{\mathcal{F}(B)}$ for all $A, B \in \mathcal{G}_{\mathrm{sym}}$ \cite{Esmi2018}.
\end{proof}

Following Section~3 of Esmi et al.\ \cite{Esmi2018}, a trapezoidal fuzzy number $A = (a; b; c; d)$ with $a \le b \le c \le d$ has membership function
\begin{equation}
\mu_A(x) = \begin{cases}
1, & x \in [b, c], \\
0, & x \le a \text{ or } x \ge d, \\
\dfrac{x - a}{b - a} \wedge \dfrac{d - x}{d - c}, & \text{otherwise}.
\end{cases}
\end{equation}
By Lemma~\ref{lem:symmetry}, $A = (a; b; c; d)$ is symmetric if and only if $b - a = d - c$ \cite{Esmi2018}. Parameterizing a symmetric trapezoidal fuzzy number by its symmetry center $m \in \mathbb{R}$, core half-width $\delta \ge 0$, and flank spread $\ell > 0$, we write
\begin{equation}
A = (m - \delta - \ell;\; m - \delta;\; m + \delta;\; m + \delta + \ell),
\end{equation}
which corresponds to the extended symmetric $L$-$L$ trapezoidal representation $(m - \delta, m + \delta, \ell, \ell)_{LL}$ with $L(u) = \max(0, 1 - u)$ \cite{Dubois1980,Bede2013} and $\alpha$-levels
\begin{equation}
[A]_\alpha = \left[m - \delta - \ell(1 - \alpha),\; m + \delta + \ell(1 - \alpha)\right], \qquad \alpha \in [0, 1].
\end{equation}
Because symmetric trapezoidal fuzzy numbers depend on two independent scale parameters ($\delta$ and $\ell$), the family of all symmetric trapezoids partitions into equivalence classes indexed by the core-to-flank ratio $\rho = \delta / \ell \ge 0$.

\begin{theorem}[Characterization of $\mathbb{R}_{\mathcal{F}(A)}$ for Symmetric Trapezoidals]
\label{prop:trap_char}
Let $A = (m_A - \delta_A - \ell_A;\; m_A - \delta_A;\; m_A + \delta_A;\; m_A + \delta_A + \ell_A)$ be a symmetric trapezoidal fuzzy number with $\delta_A \ge 0$ and $\ell_A > 0$, and let $\rho = \delta_A / \ell_A$. Denote by $\mathcal{T}_{\mathrm{sym}}(\rho)$ the family of all symmetric trapezoidal fuzzy numbers with flank spread $\ell > 0$ and core half-width $\delta = \rho \ell$. Then:
\begin{equation}
\mathbb{R}_{\mathcal{F}(A)} = \mathcal{T}_{\mathrm{sym}}(\rho) \cup \mathbb{R},
\end{equation}
and $\mathbb{R}_{\mathcal{F}(A)} = \mathbb{R}_{\mathcal{F}(B)}$ for all $A, B \in \mathcal{T}_{\mathrm{sym}}(\rho)$.
\end{theorem}

\begin{proof}
Let $B \in \mathcal{T}_{\mathrm{sym}}(\rho)$ be arbitrary with center $m_B \in \mathbb{R}$, flank spread $\ell_B > 0$, and core half-width $\delta_B = \rho \ell_B$, so that
\begin{equation}
[B]_\alpha = \left[m_B - \rho\ell_B - \ell_B(1 - \alpha),\; m_B + \rho\ell_B + \ell_B(1 - \alpha)\right], \qquad \alpha \in [0, 1].
\end{equation}
Setting $q = \ell_B / \ell_A > 0$ and $r = m_B - (\ell_B / \ell_A)m_A$, and noting that $\delta_A = \rho\ell_A$, we obtain for all $\alpha \in [0, 1]$:
\begin{align}
[\Psi_A(q, r)]_\alpha &= q[A]_\alpha + r \nonumber \\
&= \frac{\ell_B}{\ell_A}\left[m_A - \rho\ell_A - \ell_A(1 - \alpha),\; m_A + \rho\ell_A + \ell_A(1 - \alpha)\right] + \left(m_B - \frac{\ell_B}{\ell_A}m_A\right) \nonumber \\
&= \left[m_B - \rho\ell_B - \ell_B(1 - \alpha),\; m_B + \rho\ell_B + \ell_B(1 - \alpha)\right] = [B]_\alpha.
\end{align}
Thus $B = \Psi_A(q, r) \in \mathbb{R}_{\mathcal{F}(A)}$, and adjoining crisp real numbers via $\chi_{\{r\}} = \Psi_A(0, r)$ \cite{Esmi2018,Pedro2019} yields $\mathcal{T}_{\mathrm{sym}}(\rho) \cup \mathbb{R} \subseteq \mathbb{R}_{\mathcal{F}(A)}$.

Conversely, let $C = \Psi_A(q, r) \in \mathbb{R}_{\mathcal{F}(A)}$ for $(q, r) \in \mathbb{R}^2$. By Corollary~\ref{cor:linear_preservation}(3), if $q = 0$, then $C = \chi_{\{r\}} \in \mathbb{R}$ \cite{Pedro2019}. If $q \neq 0$, then $C$ is an extended symmetric $L$-$L$ trapezoidal fuzzy number with center $m_C = q m_A + r$, flank spread $\ell_C = |q|\ell_A > 0$, and core half-width $\delta_C = |q|\delta_A = \rho\ell_C$, so $C \in \mathcal{T}_{\mathrm{sym}}(\rho)$. Therefore $\mathbb{R}_{\mathcal{F}(A)} = \mathcal{T}_{\mathrm{sym}}(\rho) \cup \mathbb{R}$, and Proposition~\ref{prop:invariance} yields $\mathbb{R}_{\mathcal{F}(A)} = \mathbb{R}_{\mathcal{F}(B)}$ for all $A, B \in \mathcal{T}_{\mathrm{sym}}(\rho)$ \cite{Esmi2018}.
\end{proof}

\begin{corollary}[Characterization for Symmetric Triangular Fuzzy Numbers]
\label{cor:tri_char}
Let $\mathcal{T}_{\mathrm{tri,sym}}$ denote the family of all non-degenerate symmetric triangular fuzzy numbers $A = (m_A - \ell_A;\; m_A;\; m_A + \ell_A)$ with $m_A \in \mathbb{R}$ and $\ell_A > 0$. Then for every $A \in \mathcal{T}_{\mathrm{tri,sym}}$,
\begin{equation}
\mathbb{R}_{\mathcal{F}(A)} = \mathcal{T}_{\mathrm{tri,sym}} \cup \mathbb{R},
\end{equation}
and all non-degenerate symmetric triangular fuzzy numbers are mutually linearly correlated.
\end{corollary}

\begin{proof}
Viewing each triangular fuzzy number $(a; b; d)$ as the zero-core trapezoidal fuzzy number $(a; b; b; d)$ \cite{Esmi2018}, we have $\delta_A = 0$ and $\rho = \delta_A / \ell_A = 0$. Setting $\rho = 0$ in Theorem \ref{prop:trap_char} gives $\mathcal{T}_{\mathrm{sym}}(0) = \mathcal{T}_{\mathrm{tri,sym}}$, whence $\mathbb{R}_{\mathcal{F}(A)} = \mathcal{T}_{\mathrm{tri,sym}} \cup \mathbb{R}$.
\end{proof}

\section{Necessary and Sufficient Conditions in the Asymmetric Case}

When the generator $A \in \mathbb{R}_{\mathcal{F}} \setminus \mathbb{R}$ is non-symmetric, Esmi et al.\ \cite{Esmi2018} guarantees that $\Psi_A : \mathbb{R}^2 \to \mathbb{R}_{\mathcal{F}(A)}$ is a bijection endowed with a two-dimensional Banach space structure isomorphic to $\mathbb{R}^2$. Unlike the symmetric families in Section~4, two non-degenerate asymmetric fuzzy numbers of the same parametric type are linearly correlated if and only if their dimensionless skewness and shape invariants coincide either directly ($q > 0$) or under mirror reflection ($q < 0$).

\begin{theorem}[Linear Correlation for Asymmetric Gaussian Fuzzy Numbers]
\label{prop:asym_gauss}
Let $A = (m_A; \sigma_{L,A}, \sigma_{R,A})$ and $B = (m_B; \sigma_{L,B}, \sigma_{R,B})$ be non-degenerate asymmetric Gaussian fuzzy numbers with centers $m_A, m_B \in \mathbb{R}$ and strictly positive spreads $\sigma_{L,A}, \sigma_{R,A}, \sigma_{L,B}, \sigma_{R,B} > 0$, whose membership functions are given by
\begin{equation}
\mu_K(x) = \begin{cases}
\exp\!\left(-\dfrac{(x - m_K)^2}{2\sigma_{L,K}^2}\right), & x \le m_K, \\[10pt]
\exp\!\left(-\dfrac{(x - m_K)^2}{2\sigma_{R,K}^2}\right), & x > m_K,
\end{cases} \qquad K \in \{A, B\}.
\end{equation}
Then $A$ and $B$ are linearly correlated (with $B = \Psi_A(q, r)$ for a unique $(q, r) \in (\mathbb{R} \setminus \{0\}) \times \mathbb{R}$) if and only if
\begin{equation}
\frac{\sigma_{R,B}}{\sigma_{L,B}} = \frac{\sigma_{R,A}}{\sigma_{L,A}} \qquad \text{or} \qquad \frac{\sigma_{R,B}}{\sigma_{L,B}} = \frac{\sigma_{L,A}}{\sigma_{R,A}}.
\end{equation}
\end{theorem}

\begin{proof}
For $K \in \{A, B\}$ and $\alpha \in (0, 1]$, the $\alpha$-levels of $K$ are given by
\begin{equation}
[K]_\alpha = \left[m_K - \sigma_{L,K}\sqrt{-2\ln\alpha},\; m_K + \sigma_{R,K}\sqrt{-2\ln\alpha}\right],
\end{equation}
so $K$ is a strict $L$-$L$ fuzzy number with $L(u) = e^{-u^2/2}$, left spread $\sigma_{L,K}$, and right spread $\sigma_{R,K}$.

($\Rightarrow$) Suppose $B = \Psi_A(q, r)$ with $q \neq 0$. By Corollary~\ref{cor:linear_preservation}(2), if $q > 0$, then $\sigma_{L,B} = q\sigma_{L,A}$ and $\sigma_{R,B} = q\sigma_{R,A}$, giving $\sigma_{R,B} / \sigma_{L,B} = \sigma_{R,A} / \sigma_{L,A}$. If $q < 0$, then $\sigma_{L,B} = |q|\sigma_{R,A}$ and $\sigma_{R,B} = |q|\sigma_{L,A}$, yielding $\sigma_{R,B} / \sigma_{L,B} = \sigma_{L,A} / \sigma_{R,A}$.

($\Leftarrow$) Conversely, if $\sigma_{R,B} / \sigma_{L,B} = \sigma_{R,A} / \sigma_{L,A}$, set $q = \sigma_{L,B} / \sigma_{L,A} = \sigma_{R,B} / \sigma_{R,A} > 0$ and $r = m_B - q m_A$; Corollary~\ref{cor:linear_preservation}(2) gives $[\Psi_A(q, r)]_\alpha = [B]_\alpha$ for all $\alpha \in (0, 1]$. Similarly, if $\sigma_{R,B} / \sigma_{L,B} = \sigma_{L,A} / \sigma_{R,A}$, setting $q = -\sigma_{L,B} / \sigma_{R,A} = -\sigma_{R,B} / \sigma_{L,A} < 0$ and $r = m_B - q m_A$ yields $[\Psi_A(q, r)]_\alpha = [B]_\alpha$ for all $\alpha \in (0, 1]$. When $\sigma_{L,A} \neq \sigma_{R,A}$, uniqueness of $(q, r)$ follows directly from Theorem~3.6 of Esmi et al.\ \cite{Esmi2018}.
\end{proof}

\begin{theorem}[Linear Correlation for Asymmetric Trapezoidal Fuzzy Numbers]
\label{prop:asym_trap}
Let $A = (a_A; b_A; c_A; d_A)$ and $B = (a_B; b_B; c_B; d_B)$ be non-degenerate trapezoidal fuzzy numbers ($a_K < d_K$ for $K \in \{A, B\}$). Define for $K \in \{A, B\}$ the left flank spread $\ell_{L,K} = b_K - a_K \ge 0$, the core width $w_K = c_K - b_K \ge 0$, and the right flank spread $\ell_{R,K} = d_K - c_K \ge 0$. Then $A$ and $B$ are linearly correlated ($B = \Psi_A(q, r)$ with $q \neq 0$) if and only if there exists a scalar $\lambda > 0$ such that
\begin{equation}
(\ell_{L,B}, w_B, \ell_{R,B}) = \lambda(\ell_{L,A}, w_A, \ell_{R,A}) \qquad \text{or} \qquad (\ell_{L,B}, w_B, \ell_{R,B}) = \lambda(\ell_{R,A}, w_A, \ell_{L,A}).
\end{equation}
In particular, when $\ell_{L,A}, \ell_{R,A} > 0$, $A$ and $B$ are linearly correlated with $q \neq 0$ if and only if
\begin{equation}
\frac{w_B}{d_B - a_B} = \frac{w_A}{d_A - a_A} \qquad \text{and} \qquad \frac{\ell_{R,B}}{\ell_{L,B}} \in \left\{\frac{\ell_{R,A}}{\ell_{L,A}},\; \frac{\ell_{L,A}}{\ell_{R,A}}\right\}.
\end{equation}
\end{theorem}

\begin{proof}
For $K \in \{A, B\}$ and $\alpha \in [0, 1]$, the $\alpha$-levels of $K = (a_K; b_K; c_K; d_K)$ are given by \cite{Esmi2018}
\begin{equation}
[K]_\alpha = \left[a_K + \ell_{L,K}\alpha,\; d_K - \ell_{R,K}\alpha\right],
\end{equation}
with $[K]_1 = [b_K, c_K]$ of width $w_K = c_K - b_K$ and $[K]_0 = [a_K, d_K]$ of total support width $d_K - a_K = \ell_{L,K} + w_K + \ell_{R,K}$.

($\Rightarrow$) Suppose $[B]_\alpha = q[A]_\alpha + r$ for all $\alpha \in [0, 1]$ with $q \neq 0$, and set $\lambda = |q| > 0$.
\begin{itemize}
    \item If $q > 0$, then $a_B + \ell_{L,B}\alpha = q a_A + r + q\ell_{L,A}\alpha$ and $d_B - \ell_{R,B}\alpha = q d_A + r - q\ell_{R,A}\alpha$ for all $\alpha \in [0, 1]$. Equating coefficients of $\alpha$ gives $\ell_{L,B} = \lambda\ell_{L,A}$ and $\ell_{R,B} = \lambda\ell_{R,A}$, while evaluating at $\alpha = 1$ gives $w_B = c_B - b_B = q(c_A - b_A) = \lambda w_A$.
    \item If $q < 0$, endpoint reversal yields $a_B + \ell_{L,B}\alpha = q d_A + r - q\ell_{R,A}\alpha = (q d_A + r) + \lambda\ell_{R,A}\alpha$ and $d_B - \ell_{R,B}\alpha = (q a_A + r) - \lambda\ell_{L,A}\alpha$ for all $\alpha \in [0, 1]$, giving $\ell_{L,B} = \lambda\ell_{R,A}$, $\ell_{R,B} = \lambda\ell_{L,A}$, and $w_B = |q|w_A = \lambda w_A$.
\end{itemize}

($\Leftarrow$) Conversely, if $(\ell_{L,B}, w_B, \ell_{R,B}) = \lambda(\ell_{L,A}, w_A, \ell_{R,A})$ with $\lambda > 0$, set $q = \lambda > 0$ and $r = a_B - q a_A$. Since $d_B = a_B + \ell_{L,B} + w_B + \ell_{R,B} = q a_A + r + q(d_A - a_A) = q d_A + r$, it follows that $[B]_\alpha = [q a_A + r + q\ell_{L,A}\alpha,\; q d_A + r - q\ell_{R,A}\alpha] = q[A]_\alpha + r$ for all $\alpha \in [0, 1]$. Similarly, if $(\ell_{L,B}, w_B, \ell_{R,B}) = \lambda(\ell_{R,A}, w_A, \ell_{L,A})$, setting $q = -\lambda < 0$ and $r = a_B - q d_A$ gives $d_B = a_B + \lambda(\ell_{R,A} + w_A + \ell_{L,A}) = q d_A + r - q(d_A - a_A) = q a_A + r$, whence $[B]_\alpha = [q d_A + r - q\ell_{R,A}\alpha,\; q a_A + r + q\ell_{L,A}\alpha] = q[A]_\alpha + r$ for all $\alpha \in [0, 1]$. Finally, when $\ell_{L,A}, \ell_{R,A} > 0$, the vector proportionality condition is algebraically equivalent to invariance of the core-to-support proportion $w_K / (d_K - a_K)$ together with invariance of $\ell_{R,K} / \ell_{L,K}$ up to inversion.
\end{proof}

\begin{corollary}[Linear Correlation for Asymmetric Triangular Fuzzy Numbers]
\label{cor:asym_tri}
Let $A = (a_A; b_A; d_A)$ and $B = (a_B; b_B; d_B)$ be non-degenerate triangular fuzzy numbers with left spreads $\ell_{L,K} = b_K - a_K \ge 0$ and right spreads $\ell_{R,K} = d_K - b_K \ge 0$ ($K \in \{A, B\}$). Then $A$ and $B$ are linearly correlated ($B = \Psi_A(q, r)$ with $q \neq 0$) if and only if there exists $\lambda > 0$ such that
\begin{equation}
(\ell_{L,B}, \ell_{R,B}) = \lambda(\ell_{L,A}, \ell_{R,A}) \qquad \text{or} \qquad (\ell_{L,B}, \ell_{R,B}) = \lambda(\ell_{R,A}, \ell_{L,A}).
\end{equation}
In particular, when $\ell_{L,A}, \ell_{R,A}, \ell_{L,B}, \ell_{R,B} > 0$, $A$ and $B$ are linearly correlated with $q \neq 0$ if and only if
\begin{equation}
\frac{\ell_{R,B}}{\ell_{L,B}} = \frac{\ell_{R,A}}{\ell_{L,A}} \qquad \text{or} \qquad \frac{\ell_{R,B}}{\ell_{L,B}} = \frac{\ell_{L,A}}{\ell_{R,A}}.
\end{equation}
\end{corollary}

\begin{proof}
Viewing each triangular fuzzy number $(a_K; b_K; d_K)$ as the zero-core trapezoidal fuzzy number $(a_K; b_K; b_K; d_K)$ with $w_K = 0$ \cite{Esmi2018}, the result follows immediately from Proposition~\ref{prop:asym_trap}.
\end{proof}

\begin{remark}
\label{rem:non_symmetric}
As shown in Theorem~3.6 of Esmi et al.\ \cite{Esmi2018}, the operator $\Psi_A$ is injective if and only if $A$ is non-symmetric, in which case $(\mathbb{R}_{\mathcal{F}(A)}, +_A, \cdot)$ is a two-dimensional Banach space isomorphic to $\mathbb{R}^2$ (covering Theorems~\ref{prop:asym_gauss}--\ref{prop:asym_trap} and Corollary~\ref{cor:asym_tri} whenever the generator $A$ is asymmetric). Conversely, when $A$ is symmetric (as in Propositions~\ref{prop:gauss_char}--\ref{prop:trap_char} and Corollary~\ref{cor:tri_char}), $\Psi_A$ is two-to-one on $(\mathbb{R} \setminus \{0\}) \times \mathbb{R}$ via $\Psi_A(q, r) = \Psi_A(-q, 2qm_A + r)$, and $\mathbb{R}_{\mathcal{F}(A)}$ is equipped with the Fr\'{e}chet calculus developed in Section~3.2 of \cite{Esmi2018}.
\end{remark}

\section{Conclusion}

In this paper, we unified the level-set and joint possibility framework of $f$-correlated fuzzy numbers with the algebraic characterization of $A$-linearly correlated fuzzy number spaces $\mathbb{R}_{\mathcal{F}(A)}$. We proved that continuous, strictly monotonic $f$-correlation preserves the continuous-endpoint class $\mathcal{F}_C(\mathbb{R})$ and transforms $L$-$R$ fuzzy numbers into generalized $L_f$-$R_f$ (or $R_f$-$L_f$) fuzzy numbers, whereas nonlinear maps such as hyperbolic correlation deform piecewise-linear trapezoidal flanks into rational curves. By contrast, linear correlation ($f(x) = qx + r$) preserves the exact reference shape functions $(L, R)$ up to orientation reversal, yielding complete set-theoretic characterizations of $\mathbb{R}_{\mathcal{F}(A)}$ for symmetric Gaussian, trapezoidal, and triangular fuzzy numbers as well as sharp shape-ratio criteria in the asymmetric setting. These structural invariants provide theoretical guarantees for preserving parametric membership families along trajectories of interactive fuzzy systems.

\section*{Acknowledgments}
This study was financed in part by the Coordenação de Aperfeiçoamento de Pessoal de Nível Superior – Brasil (CAPES) – Finance Code 001.

\section*{AI Use Statement}

During the preparation of this manuscript, the authors used Google Gemini 4 Hydrogen to assist with structural organization, algebraic verification of proofs, bibliographic formatting, and \LaTeX{} typesetting, and DeepSeek V3 for text revision. All mathematical definitions, theorems, derivations, and conclusions were independently verified, edited, and approved by the authors, who take full responsibility for the scientific accuracy and integrity of the content.

\bibliographystyle{plain}
\bibliography{references}

\end{document}